\documentstyle[12pt]{article}
\title{The height problem in  first passage percolation
} 
\author{Yu Zhang
\\
Department of Mathematics, University of Colorado}
\date{}
\begin{document}
\baselineskip .20in
\maketitle

\begin{abstract}
We consider  the first passage percolation model in ${\bf Z}^2$ with a  distribution $F$ for $ 0<F(0) < p_c$. In this paper,  we  solve the height problem.
\end{abstract}

\section {Introduction of the model and results.}

We consider the first passage percolation model on  the $({\bf Z}^2, {\bf E}^2)$ lattice, a graph, with the vertices in ${\bf Z}^2$, and  with the edges in ${\bf E}^2$ connecting each pair of vertices 
one unit apart. 
We assign independently to each edge a non-negative {\em passage time} $t(e)$ with a distribution $F$ in $[0, \infty)$. More formally, we take $\Omega=[0, \infty)^{{\bf E}^2}$ as a sample space, whose points 
are called {\em configurations}.
Let $P$ be the corresponding product measure on $\Omega$. The
expectation and variance with respect to $P$ are denoted by $E(\cdot)$
and $\sigma^2(\cdot)$. 
For any two vertices ${\bf u}$ and ${\bf v}$, 
a path $\gamma$ from ${\bf u}$ to ${\bf v}$ is an alternating sequence 
$({\bf v}_0, e_1, {\bf v}_1,... ,{\bf v}_{i-1}, e_{i}, {\bf v}_{i},... ,{\bf v}_{n-1},e_n, {\bf v}_n)$ of vertices ${\bf v}_i$ in ${\bf Z}^2$ and 
 edges $e_i$ between ${\bf v}_{i-1}$ and ${\bf v}_{i}$ in ${\bf E}^2$ with ${\bf v}_0={\bf u}$ and $ {\bf v}_n={\bf v}$. 
A path is called {\em disjoint} if ${\bf v}_i\not = {\bf v}_j$ for $i\neq j$. 
In this paper, we always consider a path to be disjoint.
Given such a path $\gamma$, we define its passage time  as 
$$T(\gamma )= \sum_{e_i\in \gamma} t(e_i).\eqno{}$$
For any two sets $A$ and $B$, we define the passage time from $A$ to $B$
as
$$T(A,B)=\inf \{ T(\gamma): \gamma \mbox{ is a path from $A$ to $B$}\},$$
where the infimum is taken over all possible finite paths from some vertex in $A$ to some vertex in $B$.
A path $\gamma$ from $A$ to $B$ with $T(\gamma)=T(A, B)$ is called an {\em optimal path} of 
$T(A, B)$.  

The existence of such an optimal path has been proven (see Kesten (1986)).
If $t(e)=0$, the edge is called a {\em zero edge} or an {\em open edge}; otherwise, it is called a {\em closed edge}.
We also want to point out that the optimal path may  not be  unique.
If all edges in a path are in passage time zero,  the path is called a {\em zero path} or
an {\em open path}.
If we focus on a special configuration $\omega$, we may write $T(A, B)(\omega)$
instead of $T(A, B)$. 
When $A=\{\bf u\}$ and $B=\{\bf v\}$ are single vertex sets, $T({\bf u}, {\bf v})$ is the passage time
from ${\bf u}$ to ${\bf v}$. We may extend the passage time over ${\bf R}^2$.
More precisely,
if ${\bf u}$ and ${\bf v}$ are in ${\bf R}^2$, we define $T({\bf u}, {\bf v})=T({\bf u}', {\bf v}')$, where
${\bf u}'$ (resp., ${\bf v}'$) is the nearest neighbor of ${\bf u}$ (resp., ${\bf v}$) in
${\bf Z}^2$. Possible indetermination can be eliminated by choosing an
order on the vertices of ${\bf Z}^2$ and taking the smallest nearest
neighbor for this order.  In this paper, for any ${\bf x}, {\bf y}\in {\bf R}^2$,
$\|{\bf x}\|$ is denoted by the Euclidean norm and
$d({\bf x}, {\bf y})=\|{\bf x}-{\bf y}\|$ is the {\em distance} between ${\bf x}$ and ${\bf y}$.
For any two sets ${ A}$ and ${ B}$ of ${\bf R}^2$, 
$$d({A}, {B})=\min\{ d({\bf x}, {\bf y}): {\bf x}\in {A}\mbox{ and } {\bf y}\in {B}\}$$
and is denoted as the distance between ${A}$ and ${B}$.
For a vertex set $A$, we denote by $|A|$  the number of vertices in $A$. Similarly, if $A$ is an edge set, we also denote by $|A|$  the number of edges in $A$.

Given ${\bf 0}=(0, 0)$ and ${\bf x}=(x_1, x_2)\in {\bf R}^2$, if $Et(e) < \infty$, by Kingman's  sub-additive  ergodic theorem,  it is well known that
$$
\lim_{n\rightarrow  \infty}{1\over n} T({\bf 0}, n{\bf x}) = \inf_{n} {1\over n} E T({\bf 0}, n{\bf x})=\lim_{n\rightarrow \infty}
{1\over n} E T({\bf 0}, n {\bf x})=\mu_F({\bf x}) \mbox{ a.s. and in } L_1.\eqno{(1.1)}
$$
It is also known (see Kesten (1986))  that
$$\mu_F({\bf x}) \mbox{ is continuous in } {\bf x}\mbox{ and } \mu_F({\bf x})>0\mbox{ for ${\bf x}\neq 0$ iff }F(0)  < 1/2.\eqno{(1.2)}$$
In particular, Hammersley and Welsh (1965), in their pioneering  paper, investigated a special case ${\bf x}=(1,0)$:
 $$a_{0, n}=T({\bf 0}, n{\bf x})=T({\bf 0}, (n,0)).$$
They showed  that
$$\lim_{n\rightarrow \infty}  a_{0, n}/n  =\mu_F((1,0))\mbox{ a.s. and in }L_1.\eqno{(1.3)}$$
For simplicity's sake, we denote 
$$\mu_F(1,0)=\mu_F.\eqno{}$$

Now let  ${\bf L}[n]$ be the vertices in the  vertical line  passing through $(n, 0)$.  We look again at Hammersley and Welsh (1965)  and their introduction of  point-line passage time. 
Let
$$b_{0,  n}=\inf\{ T(\gamma): \gamma  \mbox{ is a path from the origin to ${\bf L}[n]$}\}.$$
   They proved in their paper
   $$\lim_{n\rightarrow \infty} b_{{\bf 0}, n}/n=\mu_F \mbox{ a.s. and in }L_1.\eqno{(1.4)}$$
Since we need the rate  of convergence in (1.1), we assume in this paper  that $t(e)$ is not a constant and satisfies the following:
$$\int_0^\infty e^{\iota x} dF(x) < \infty \mbox{ for some } \iota >0.\eqno{(1.5)}$$
When $F(0) < 1/2$, the map ${\bf x} \rightarrow \mu_F({\bf x})$ induces  a norm on
${\bf R}^2$. The unit radius ball for this norm is denoted by 
$${\bf B}:={\bf B}(F)=\{{\bf u}: \mu_F({\bf u})\leq 1\}$$
and is called the {\em asymptotic shape}. The boundary of ${\bf B}$ is
$$\partial {\bf B}:= \{ {\bf u} \in {\bf R}^2: \mu({\bf u})=1\}.$$
By (1.2), if $F(0) < 1/2$, ${\bf B}$ is a compact convex
deterministic set.
Define for all $t> 0$ a {\em random shape} 
$${\bf B}(t):= \{{\bf u}\in {\bf Z}^d, \ T( {\bf 0}, {\bf u}) \leq t\}.$$
The shape theorem (see Cox and Durrett (1981)) is a well-known result stating that for any
$\epsilon >0$,
$$t(1-\epsilon){\bf B}  \subset {{\bf B}(t) } \subset t(1+\epsilon ){\bf B}
\mbox{ eventually w.p.1.}\eqno{(1.6)}$$

It is easy to show (see Fig. 6.1 in Kesten (1986)) that
\begin{eqnarray*}
&&{\bf B}\mbox{ lies  between the cube $[-1/\mu_F, 1/\mu_F]^2$ and 
 the diamond with
the four corners }\\
&&\mbox{$( \pm 1/\mu_F, 0)$ and $(0, \pm 1/\mu_F)$}.\hskip 11.5cm {(1.7)}
\end{eqnarray*}
By using a rate estimate (see  Chow and Zhang (2003)), if $F(0) < p_c$ and (1.5) holds, for any $\epsilon >0$, then
$$P(t(1-\epsilon){\bf B}  \subset {{\bf B}(t) } \subset t(1+ \epsilon ){\bf B})\geq 1- \exp(-O(\epsilon t)).\eqno{(1.8)}$$
By Zhang's (2010) Theorem 3, if $F(0) < 1/2$ and (1.5) holds, then
 $$P(|b_{{\bf 0}, n}-n \mu_F|\geq \epsilon n )\leq \exp(-O( \epsilon n)).\eqno{(1.9)}$$
 Kesten (1996) also showed the following concentration inequality. If $F(0) < 1/2$ and (1.5) holds, then for any $z>0$, there exist $C_i=C_i(F)$ for $i=1,2$ independent of   ${\bf x}$ such that
 $$ P(|T({\bf 0}, {\bf x})-ET({\bf 0}, {\bf x})|\geq z \|{\bf x}\|^{1/2})\leq C_1\exp(-C_2 z^2 ).\eqno{(1.10)}$$
In this paper, $C_i$ denotes a constant with $0< C_i < \infty$ whose precise value is of no importance; its value may change from appearance to appearance, but $C_i$ will always be independent  of  $n$ and $t$, $k$, and $m$, although it  may depend on $F$.
For simplicity's sake, we sometimes use $O(n)$ for $C_1 n \leq O(n) \leq C_2 n$ if we do not need the precise value of $C_i$. 
Regarding  the length of optimal paths, if $F(0) < 1/2$, it is well known (see Prop. (5.8) in Kesten (1986)) that  there exists $\lambda=\lambda(F)$ and $C_i=C_i(F)$ for $i=1,2$ such that
$$P(\exists \, \mbox{ an optimal path $\gamma_n({\bf x})$ with }|\gamma_n({\bf x})|\geq \lambda n)\leq C_1\exp(-C_2 n).\eqno{(1.11)}$$
 
We know that optimal paths of $a_{0, n}$ exist, but they are not unique.  When $F(0) < 1/2$,  each optimal  path of $a_{0, n}$ is finite.
Similarly, optimal paths of $b_{0, n}$ exist, but they also are not unique.  When $F(0) < 1/2$,  each optimal  path of $b_{0, n}$ is finite.
 We denote by   $A_n$   and $B_n$ the unions of all optimal paths with passage times $a_{0, n}$ and  $b_{0, n}$, respectively.
 Moreover, $\bar{A}_n$ and $\bar{B}_n$ are  the set of all the optimal paths with passage times $a_{0, n}$ and  $b_{0, n}$, respectively.
 $|\bar{A}_n|$ and $|\bar{B}_n|$ are the numbers of paths in $\bar{A}_n$ and $\bar{B}_n$, respectively.
We denote the {\em height } of the optimal paths of $a_{0, n}$ or $b_{0, n}$  by 
$$h_n^\alpha=\max\{ d({\bf u},  \mbox{the $X$-axis}): {\bf u}\in \alpha_n\}\mbox{ for }\alpha=A \mbox{ or }\alpha=B.$$
 It is widely believed  (see Hammersley and  Welsh  (1965); Smythe and Wierman  (1978); and 
Kesten (1986))  that
$$\lim _{n\rightarrow \infty} {h_n^\alpha\over n} =0 \mbox{ in some sense for }\alpha=A \mbox{ or }\alpha=B \mbox{ if } F(0) < 1/2.\eqno{(1.12)}$$
The conjecture of  (1.12) is called the {\em height problem}.   Note that  (1.12) does not hold when $F(0) \geq 1/2$.
In this paper, we answer this conjecture affirmatively. \\
 
 {\bf Theorem.} {\em If  $F$ is a  distribution with $0< F(0) <1/2$ and (1.5) holds, then 
$$\lim _{n\rightarrow \infty} {h_n^\alpha \over n} =0 \mbox{ in probability and in $L_1$ for }\alpha=A \mbox{ or }\alpha=B.$$}

 {\bf Remarks 1-4.}  1. In this paper,  we only show the theorem for $a_{0, n}$ since the geometry is much easier to handle. In other words, we solve the height conjecture in the horizontal direction. 
 But the proof for the theorem might be carried over to solve the height conjecture in any direction.
  
  2.  A distribution F is said to be {\em useful} (see van den Berg and Kesten (1993)) if $F(0) < 1/2$ and $F^{-}=0$ or   $F(F^{-})< \vec{p_c}$, where 
$\vec{p}_c$ is the critical probability of two-dimensionally oriented percolation and $F^-$ is the infimum of the support of $F$. The condition that $0< F(0) < 1/2$ in the theorem  is a special  case of useful distribution. 
The same proof of the  theorem might  show the theorem  when $F$ is useful and (1.5) holds.   We also want to point out that  there is a counter example (see Durrett and Liggett (1981)) that 
 the theorem does not hold  in some direction when $F$ is not useful.

3. We may define the {\em height fluctuation exponent} to be
$$ \xi(d)=\inf\left\{ \xi:  \limsup_n{E  h_n^A\over n^{\xi} }<\infty\right\}.$$
 It is believed that $\xi(2)=2/3$ if $F(0) < 1/2$. 
We proved in  the theorem  that the ratio of the height goes to zero. 
But we are unable to show that
$$\lim _{n\rightarrow \infty} {h_n^A\over n^{1-\delta}} =0 \mbox{ in probability}\eqno{}$$
for some $\delta >0$.

4. Our proof for the theorem only works for  the two-dimensional lattice. \\

\section{ The number of pivotal edges  for optimal paths.}
Among most percolation problems, people have to focus on some specific edges, called pivotal edges, which optimal paths have to use.
We first introduce the definition of pivotal edges.
 We denote by $L[i]$ the vertical line passing $(i, 0)$. We also denote by $L[i, j]$  ($L(i, j)$) the  space between $L[i]$ and $L[j]$, called   a
{\em cylinder}, including $L[i]$ and $L[j]$ (but not including $L[i]$ and $L[j]$) for $i< j$. Let  ${\bf L}[i, j]$  be the vertex set of $L[i, j]$. When $i=j$, ${\bf L}[i]$ is the vertex set of $L[i]$ defined before.
 Hammersley and Welsh (1965) introduced the cylinder  point-line passage time:
 \begin{eqnarray*}
 &&s_{0,  n}=\inf\{ T(\gamma): \gamma  \mbox{ is a path from the origin to ${\bf L}[n]$ with vertices in  ${\bf L}(0,n)$}\\
 &&\hskip 2cm \mbox{  except its initial and terminate vertices}\},
 \end{eqnarray*}
   and they proved 
   $$\lim_{n\rightarrow \infty} s_{{0}, n}/n=\mu_F \mbox{ a.s. and in }L_1.\eqno{(2.1)}$$
   By the same proof of Zhang's (2010) Theorem 3, if $F(0) < p_c$ and (1.5) holds, then
 $$P(|s_{0, n}-n \mu_F|\geq \epsilon n )\leq \exp(-O( \epsilon n)).\eqno{(2.2)}$$
For convenience,   we give a general  definition of cylinder point-line passage time starting at any vertex as follows. For $a < n$, 
   \begin{eqnarray*}
   &&s_{a, n} ((a, b))=\inf\{ T(\gamma): \gamma  \mbox{ is a path from $(a, b)$  to ${\bf L}[n]$  with vertices in  ${\bf L}(a,n)$}\\
   &&\hskip 3cm \mbox{ except its inital and terminal vertices}\}.
   \end{eqnarray*}
 We need to discuss  a few properties regarding  optimal path $\gamma_n$ of $s_{0, n}$.  It follows from Proposition (5.8) in Kesten (1986) that  there exists $\lambda=\lambda(F)$ and $C_i=C_i(F)$ for $i=1,2$ such that
 $$P(\exists\,\, \gamma_n \mbox{ with } T(\gamma_n)= s_{0, n} \mbox{ such that }|\gamma_n|\geq \lambda m \mbox{ with }m\geq n) \leq \exp(-O(m)).\eqno{(2.3)}$$
 Let $S_n$ be the union of  all the optimal paths with passage time $s_{0, n}$. We also denote by $\bar{S}_n$  the set of all optimal paths with passage time $s_{0, n}$.  Furthermore, let  $D_n$ be the intersection of all the optimal paths in $\bar{S}_n$.
 For $e\in D_n$, each optimal path in $\bar{S}_n$ has to go through $e$, so
      the edges in $D_n$ are called {\em pivotal edges}.
By using  Theorem 2 in Zhang (2006),  if $F$ is a Bernoulli distribution with $F(0) < 1/2$, then  there exists $M$ such that
$$ E|S_n|\leq M n.\eqno{(2.4)}$$
Later,  Nakajima (2019) showed that (2.4) holds for a point-point passage time and for  a useful distribution. Furthermore,  he also showed that  for a point-point passage time,
$$E|D_n|=O(n).\eqno{(2.5)}$$

Now we need to account for the number of pivotal edges in a slab.
We first investigate the behavior of optimal paths for $s_{0, n}$ in a slab.  We consider an optimal path $\gamma_{n}\in \bar{S}_n$.
  Let $0<  m_1< m_2<n$ with $m_1=O(n)$,  $m_2-m_1=\kappa n$ and $n-m_2=O(n)$ for a small, positive constant $\kappa$.
$\gamma_n$  first meets   ${\bf u}'$ at  $L[{m_1}]$ and   last meets  ${\bf u}$ at $L[{m_1}]$. Furthermore,  $\gamma_n$ continues to first meet  ${\bf v}$ at $L[m_2]$ and last meet ${\bf v}'$ at $L[m_2]$.  Finally, $\gamma_n$ meets  $L[n]$ at ${\bf z}$. We denote by $\gamma({\bf 0}, {\bf u'})$
the sub-path from ${\bf 0}$ to ${\bf u}$. Similarly, we denote the sub-paths $\gamma({\bf u}', {\bf u})$,  $\gamma({\bf u}, {\bf v})$, $\gamma({\bf v}, {\bf v}')$,   as well as  $\gamma({\bf v}', {\bf z})$ the corresponding  point-point sub-paths.   We sometimes need to  assume that ${\bf v}={\bf z}$. 
We now  show the following lemma.\\

{\bf Lemma 2.1.} {\em If  $F$ is a  distribution with $F(0) < 1/2$, $0< \kappa< 1$ is a  constant, and (1.5) holds, then for any large $n$, $0<m_1< m_2 < n$ with $m_1=O(n)$, $m_2-m_1=\kappa  n$, and $n-m_2=O(n)$, and for any  $\gamma_{n}\in \bar{S}_n$ and  for  any $0<\delta < \mu_F/2$, there exist $C_j=C_j(F, \delta)$ for $j=1,2$ such that
$$P(\forall\,\,\gamma_n\in \bar{S}_n, \kappa n (\mu_F+\delta ) \geq T(\gamma({\bf u}', {\bf v}'))\geq T(\gamma({\bf u}, {\bf v}))\geq  \kappa n (\mu_F-\delta))\geq 1-C_1\exp(- C_2\kappa n).\eqno{}$$}

{\bf Proof.}  
 We denote the three sub-paths of  $\gamma_n$ from the origin to ${\bf u}'$,   from ${\bf u}'$ to ${\bf v}'$,  and from ${\bf v}'$ to ${\bf z}$
 by
$$\gamma_n=\gamma({\bf 0}, {\bf u}')\cup \gamma({\bf u}', {\bf v}')\cup \gamma({\bf v}', {\bf z}).\eqno{}$$
Note that  these  sub-paths are also optimal paths from ${\bf 0}$ to $ {\bf u}'$, from ${\bf u}'$ to $ {\bf v}'$,  and from ${\bf v}'$ to ${\bf z}$, respectively.
Thus, for any $0< \delta \leq \mu_F/2$,   by  (2.3) and  (1.8),  there exist $C_l=C_l(F(0))$ for $l=1,2,$
\begin{eqnarray*}
&&P(\exists\,\,\gamma_n\in \bar{S}_n,(\mu_F-\delta  /8)m_1 \geq  T(\gamma({\bf 0}, {\bf u}')))\\
&\leq & P(\exists \,\,\gamma_n\in \bar{S}_n,(\mu_F-\delta /8)m_1 \geq  T(\gamma({\bf 0}, {\bf u}')),|\gamma_n|\leq \lambda n)+C_1\exp(-C_2 n)\\
&\leq &\sum_{\bf w}P(\exists \,\,\gamma_n\in \bar{S}_n,(\mu_F-\delta/8)m_1 \geq  T(\gamma({\bf 0}, {\bf w})),|\gamma_n|\leq \lambda n, {\bf u}'={\bf w})+C_1\exp(-C_2 n)\\
&\leq &( \lambda n)^d \exp(-O(\delta n))+C_1\exp(-C_2 n)\leq \exp(-O(\delta n)) ,\hskip 6.5cm {(2.6)}
\end{eqnarray*}
 where the  sum  takes all possible  ${\bf w}\in {\bf L}[m_1]$  with 
$\|{\bf w}\|\leq 2\lambda n$, and $\gamma({\bf 0}, {\bf w})$ is an optimal path from ${\bf 0}$ to ${\bf w}$.
Thus,  by (2.6),
$$P(\forall\,\,\gamma_n\in \bar{S}_n,(\mu_F-\delta/8)m_1 \leq T(\gamma({\bf 0}, {\bf u}')))\geq 1-\exp(-O(\delta n)).\eqno{(2.7)}$$
Similarly, we have 
$$ P(\forall\,\,\gamma_n\in \bar{S}_n,(\mu_F-\delta /8)  (n-m_2)\leq T( \gamma({\bf v}', {\bf z})))\geq 1-\exp(-O(\delta n)). \eqno{(2.8)}$$

On the other hand,  for ${\bf u}$ and ${\bf v}$, by  (2.2), 
$$1-e^{-O(\delta n)} \leq  P(\forall\,\,\gamma_n\in \bar{S}_n,T( \gamma_n)=T(\gamma({\bf 0}, {\bf u}'))+T(\gamma({\bf u}', {\bf v}'))+ T(\gamma({\bf v}', {\bf z}))\leq n(\mu_F+\delta/8)).\eqno{(2.9)}
$$
Substituting (2.7) and (2.8)  into (2.9),  note that $m_2-m_1=\kappa n$, so
$$ 1-\exp(-O(\delta \kappa n))\leq P(\forall\,\,\gamma_n\in \bar{S}_n,T(\gamma({\bf u}', {\bf v}'))\leq  n \kappa(\mu_F+\delta )).\eqno{(2.10)}$$
 Note that $d({\bf u}, {\bf v}) \geq n\kappa$, so  by the same proof of (2.7), 
 $$ 1-\exp(-O(\delta \kappa n))\leq P(\forall\,\,\gamma_n\in \bar{S}_n,T(\gamma({\bf u}, {\bf v}))\geq  n \kappa(\mu_F-\delta )).\eqno{(2.11)}$$
 Lemma 2.1 follows from (2.10) and (2.11). $\Box$\\

Regarding the length of $\gamma({\bf u}', {\bf v}')$, by using Proposition (5.8) in Kesten (1986) and the inequality in Lemma 2.1, we have the following lemma.\\

{\bf Lemma 2.2.} {\em If  $F$ is a  distribution with $F(0) < 1/2$,  $0< \kappa< 1$ is a  constant, and (1.5) holds, then for any large $n$, $0<m_1< m_2 < n$ with $m_1=O(n)$, $m_2-m_1=\kappa  n$, and $n-m_2=O(n)$, there exist $\lambda=\lambda(F)$ and $C_j=C_j(F)$ for $j=1,2$ such that for all $m\geq n$,
$$P(\forall\,\,\gamma_n\in\bar{S}_n, |\gamma({\bf u}', {\bf v}') |\geq \lambda \kappa m)\leq C_1\exp(- C_2\kappa m).$$}

We want to remark that the same proofs of Lemmas 2.1 and 2.2 also work on  the case  when ${\bf v}={\bf z}$.
   For each optimal path $\gamma_n\in \bar{S}_n$, we can take $m_1=n-\kappa n$ and $m_2=n$ and the corresponding ${\bf u}'$, ${\bf u}$, and ${\bf v}={\bf z}$ defined in Lemma 2.1 for a small $0 < \kappa<1$.  In this special case, Lemma 2.1 and Lemma 2.2 hold by the same proofs directly.
   Recall that $D_n$ is the number of pivotal edges  of $S_n$. We  denote by  $ D_n(\kappa)$  the  number of the  pivotal edges of $S_n$ in $L[n-m_1, n]$. 
   We select any optimal path $\gamma_n$ and consider the pivotal edges of $S_n$ in $\gamma({\bf u}, {\bf v})$. 
   We may select another optimal path $\bar{\gamma}_n$ and
    consider the pivotal edges of $S_n$ in $\bar{\gamma}(\bar{\bf u}, \bar{\bf v})$.   
    Since all optimal paths have to pass all the pivotal edges,  the pivotal edges in $\gamma({\bf u}, {\bf v})$ and the pivotal edges in $\bar{\gamma}(\bar{\bf u}, \bar{\bf v})$ are the same.     Thus, we may study the  pivotal edges of $S_n$ in $L[n-m_1, n]$ in $\gamma({\bf u}, {\bf v})$  by using a particular optimal path $\gamma_n$. 
Now we  use  Lemma 2.1,  the method of  Lemma 8 in  Nakajima (2019), and the  Peierls argument  to estimate the pivotal edges  in the sub-optimal path $\gamma({\bf u}, {\bf v})$ of  $\gamma_n\in \bar{S}_n$.\\

{\bf Lemma 2.3.} 
{\em If  $F$ is a   distribution with $0< F(0) < p_c$ and (1.5) holds, then for any small $0< \kappa< 1$, there exist constants $m$  and $M$  independent of $\kappa$ and $n$ such that
$$m\kappa  n \leq ED_n(\kappa)\leq M\kappa n.$$}

{\bf Proof.}  In the proof of Lemma 2.3, we use the same notations $\kappa $, ${\bf u}$, and  ${\bf v}= {\bf z}$ as we did in Lemma 2.1.
 Let $\tau(e)$ be an independent copy of $t(e)$ for each edge $e$. Given an  edge $f$, we set $\beta^f(e)$:
$$
\beta^f(e)=\left\{ \begin{array}{cc}
\tau(e) & \mbox{ if $e=f,$}\\
t(e) &\mbox{ if $e\neq f $.}
\end{array}
\right.
$$
For an edge $f\in {\bf L}[0, n]$, let $\bar{S}_n^f$ be the set of all optimal paths from the origin to ${\bf L}[n]$ in ${\bf L}(0, n)$ except their initial and terminate vertices with edge weights  $\{\beta^f(e)\}$ for $e\in {\bf L}[0, n]$.
In other words, $\bar{S}_n^f$ is the set of all optimal paths from the origin to ${\bf L}[n]$ with edge weights  $t(e)$ for $e\neq f$  and  edge weight $\tau(f)$ for $e=f$.
 Since $t(e)$ and $\beta^f(e)$ have the same distribution for each edge $e$, for a fixed edge $f\in {\bf L}[0, n]$,
\begin{eqnarray*}
&& P( \exists \,\,\gamma_n \in \bar{S}_n, f\in \gamma({\bf u}, {\bf v}),  t(f) > 0)P( \tau(f) =0)= P( \exists \,\,\gamma_n \in \bar{S}_n ,  f\in \gamma({\bf u}, {\bf v}),  t(f) > 0,  \tau(f) =0)\\
&&= P( \exists \,\,\gamma_n \in \bar{S}_n,  f\in \gamma({\bf u}, {\bf v}),  t(f) > 0,  \beta^f(f) =0)\leq  P( \forall  \,\,\gamma_n \in \bar{S}^f_n,  f\in \gamma({\bf u}, {\bf v}))\\
&&= P( \forall  \,\,\gamma_n \in \bar{S}_n,  f\in \gamma({\bf u}, {\bf v})).\hskip 11.2cm (2.12)
\end{eqnarray*}
Thus, the event in the right side of (2.12) implies that $f$  is a pivotal edge of $\gamma_n$.  For each $\gamma_n\in \bar{S}_n$, let 
$$N_n(\gamma_n, {\bf u}, {\bf v})= \mbox{ the number of   $f$ such that }   f\in \gamma({\bf u}, {\bf v}) \mbox{ with $t(f) >0$}.$$
By (2.12), if we denote by
$$N_n(\kappa)= \min_{\gamma_n\in \bar{S}_n} \{N_n(\gamma_n, {\bf u}, {\bf v})\},$$
then 
\begin{eqnarray*}
&& ED_n(\kappa)= \sum_{f\in {\bf L}[0, n]}  P(  \forall  \,\,\gamma_n \in \bar{S}_n,  f\in \gamma({\bf u}, {\bf v}))\\
&\geq & P( \tau(f)=0)  \sum_{f\in {\bf L}[0, n]}  P(  \exists \,\,\gamma_n \in \bar{S}_n,  f\in \gamma({\bf u}, {\bf v}), t(f) >0)\geq  F(0)  EN_n(\kappa).\hskip 3cm (2.13)
\end{eqnarray*}

Now we estimate $EN_n(\kappa)$. Note that for a small $\delta >0$,
$$\{N_n(\kappa) \leq  \delta \kappa n\}\leq \{\exists\,\, \gamma_n\in \bar {S}_n \mbox{ such that } N_n(\gamma_n, {\bf u}, {\bf v}  )\leq \delta \kappa n\}.\eqno{(2.14)}$$
Note that the number of  optimal paths in $\bar { S}_n$ is finite, so on $\{N_n(\gamma_n, {\bf u}, {\bf v}  )\leq \delta \kappa n\}$, there is a $\gamma_n\in \bar{S}_n$ such that
$$N_n(\gamma_n, {\bf u}, {\bf v})\leq  \delta \kappa n.$$
If there are many optimal paths satisfying above inequality, we simply select one in a unique way.
For ${\bf z}=(z_1, z_2)\in {\bf Z}^2$ and  $M$ large independent of $n$, we denote $M$-squares by
$$\{B_{M} ({\bf z})\}=\{[ M z_1, M z_1+M]\times [ M z_2, M z_2+M]\}.$$
 For the selected optimal path $\gamma_n\in \bar{S}_n$, if $B_M({\bf z})\cap \gamma({\bf u}, {\bf v})\neq \emptyset$, then $B_M({\bf z})$ is called a {\em good   square}.  Given a good square, we consider its $3M$-squares, called  3$M$-good squares,  each containing   $B_M({\bf z})$ as its center square and eight $M$-neighbor squares of $B_M({\bf z})$. Since $\gamma({\bf u}, {\bf v})$ crosses the annulus between an $M$-good square and the  boundary of its $3M$-square, 
$$\mbox{ any $3M$-good square contains at least $M$ vertices of $\gamma({\bf u}, {\bf v})$}.\eqno{(2.15)}$$
We say two squares are adjacent if they have a common vertex. 
Since $ \gamma({\bf u}, {\bf v})$ is a path, it is easy to verify that all $M$-good squares are adjacent. Let $\gamma({\bf u}, {\bf v}, M)$ be 
these $M$-good squares, and let $|\gamma({\bf u}, {\bf v}, M)|$ be the number of $M$-squares in $\gamma({\bf u}, {\bf v}, M)$.  Note that there are nine many $M$-squares in each $3M$-square, so  for each $M$-good  square in $\gamma({\bf u}, {\bf v}, M)$, its good $3M$-square  corresponds to $M$ disjoint vertices in $\gamma({\bf u}, {\bf v})$.  Thus,
$$|\gamma({\bf u}, {\bf v}, M)|\geq 9|\gamma({\bf u}, {\bf v})|/M\geq 9\kappa n/M.\eqno{(2.16)}$$
  For the selected  $\gamma_n\in \bar{S}_n$, there are  many choices for different locations of $\gamma({\bf u}, {\bf v})$. On $|\gamma_n|\leq \lambda n$,  there are at most $(\lambda n)^2$ many  choices for the locations for vertex ${\bf u}$.   If ${\bf u}$ is fixed with ${\bf u}\in B_M({\bf z})$ for some ${\bf z}$, and if $|\gamma({\bf u}, {\bf v}, M)|=k$, by (4.24) in Grimmett (1999), we note that $\gamma({\bf u}, {\bf v}, M)$ is adjacent, so
  $$\mbox{  there are  at most   $7^{2k}$ many choices  for the locations of these $M$-squares in }\gamma({\bf u}, {\bf v}, M).\eqno{(2.17)}$$
We say a good $M$-square is {\em regular} if its $3M$-square contains an edge $e\in \gamma({\bf u}, {\bf v})$ with  $t(e) >0$; otherwise, we  say it is {\em irregular}.
Thus, if an $M$-square is irregular,  by (2.15), there is an open path in its $3M$-square with  at least $M$ vertices. Since $F(0) < 1/2$,  by Theorem 5.4 in Grimmett (1999), there exist $C_j=C_j(F(0))$ for $j=1,2$ such that
$$ P(\mbox{a fixed $M$-square is irregular}) \leq C_1\exp(-C_2M).\eqno{(2.18)}$$

With  the above observations, we use  a standard Peierls argument (see Theorem 3 in Zhang (2006) or Theorem 2.13 in van den Berg and Kesten (1993)) to estimate  the probability of the event for the number of closed edges in $\gamma({\bf u}, {\bf v})$.  By (2.13)  and (2.15), if 
$|\gamma({\bf u}, {\bf v}, M)|=k\geq 9\kappa n/ M$,  then for $0< \eta$, there exist $C_j=C_j(F)$ for $j=3, 4$,
\begin{eqnarray*}
&& P(N_n(\kappa)\leq  \eta \kappa n)\leq P(\exists \,\,\gamma_n \in \bar{S}_n, N_n(\gamma_n, {\bf u}, {\bf v})\leq  \eta \kappa n)\\
&\leq &\sum_{k\geq 9\kappa n/  M}(\lambda n)^27^{2k}\max_{\Gamma_M(k)} P( \exists \,\,\gamma_n \in \bar{S}_n,N_n(\gamma_n, {\bf u}, {\bf v})\leq  \eta \kappa n, \gamma({\bf u}, {\bf v}, M)=\Gamma_M(k)), \hskip 2cm (2.19)
\end{eqnarray*}
where  $\Gamma_M(k)$ is a fixed, adjacent $M$-square set containing ${\bf u}$ and ${\bf v}$ with $k$ many $M$-squares.
 If 
$|\gamma({\bf u}, {\bf v}, M)|=k\geq 9\kappa n/ M$, then we choose  $0 < \eta=\eta(M)\leq 9/(2M)$ small, independent of $n$ and $\kappa$,  such that 
the number of irregular $M$-squares is larger than $k/2$. 
On $\{\gamma({\bf u}, {\bf v}, M)=\Gamma_M(k)\}$,   note that there are  nine many $M$-squares in each $3M$-square, so we select disjoint, adjacent $k/18$ many  disjoint $3M$-squares such that each $3M$-square contains an irregular $M$-square selected above.  Thus, by (2.18) and (2.19),
$$P(N_n(\kappa)\leq  \eta \kappa n)\leq \sum_{k\geq \kappa n/  (9M)} (\lambda n)^27^{2k}{ k\choose k/2} (C_1 \exp(-C_2 M)) ^{k/(18)}.\eqno{(2.20)}$$
By  (2.19), note that $C_1$ and $C_2$ are independent of $M$, so by choosing $M$ large, there exist $C_j=C_j(F(0),\eta)$ for $j=3,4$ such that
$$P(N_n(\kappa)\leq  \eta \kappa n)\leq C_3 \exp(-C_4 \kappa n).\eqno{(2.21)}$$
By (2.21),  there exists $C_5$ independent  of $\kappa$ and $n$ such that 
 $$E N_n(\kappa)\geq  E (N_n(\kappa); N_n(\kappa)\geq \eta \kappa  n)\geq \eta \kappa n (1-C_3\exp(-C_4 n))\geq C_5 \kappa n   .\eqno{(2.22)}$$
 By (2.13) and (2.22), 
 $$ED_n(\kappa)\geq F(0)  C_5 \kappa n .\eqno{(2.23)}$$
On the other hand, by Lemma 2.2,
$$ ED_n(\kappa) \leq E|\max_{\gamma_n\in\bar{S}_n}\gamma({\bf u}, {\bf v})|\leq O(\kappa n).\eqno{(2.24)}$$
Lemma 2.3 follows from (2.23) and (2.24). $\Box$\\

{\bf Remark 5.}  Lemmas 2.1, 2.2, and 2.3 also hold for any $d$-dimensional  lattice.


\section{ Proof of theorem.}
In section 3, we  show two lemmas for  a point-line cylinder  height.  With these two lemmas, we show the theorem. These proofs are much easier to understand with graphs than words. We suggest that readers  use the following three graph proofs for help in the formal proofs of  the two lemmas and the theorem.\\

In their paper, Hammersley and  Welsh (1965)  introduced another height:
$$ {H}_n= \max \{|m|:  \mbox{ if the terminate vertex of $\gamma_n\in \bar{S}_n$ is $(n, m)$}\}.$$
They pointed out that the knowledge of the behavior of $ H_n$ would give information about the rate of convergence of $a_{0. n}$.
We show  the following lemma.\\

{\bf Lemma 3.1.} {\em If  $F$ is a  distribution with $0<F(0) < 1/2$ and (1.5) holds, then
$$\lim_{n\rightarrow \infty}  {{H}_n\over n}=0\mbox{ in probability}.\eqno{}$$}

{\bf Proof.}  
It is known (see Theorem 8.15 in Smythe and Wierman (1978))   that if $F(0) < 1/2$, then
$$ \limsup_{n} H_n/n \leq 1 \mbox{ a.s.}\eqno{(3.1)}$$
We suppose that Lemma 3.1 does not hold.   By  symmetry and (3.1), there is  a subsequence $\{n_s\}$ such that
$$P(1\geq H_{n_s}/ n_s \geq \delta ) \geq \delta\mbox{ for each $n_s$ and for some } \delta >0.\eqno{(3.2)}$$
For a small  $\epsilon$ with $\delta \geq 8\epsilon>0$,  by  (3.2),   there exists $l_s$ with $1\geq l_s\geq \delta$ such that for each $s$,
$$P((l_s-\epsilon) \leq  H_s/ n_s\leq (l_s+\epsilon)\mbox{ for each }s) \geq \epsilon \delta/2 .\eqno{}$$
Since $\epsilon \delta/2 $ is uniformly bounded away from zero, it is easy to find  a subsequence $\{n_{s_t}\}$ from $\{n_s\}$ independent of configurations such that 
$$P\left(\bigcap_{t=1}^\infty \{(l_{s_t}-\epsilon) \leq  H_{s_t}/ n_{s_t}\leq (l_{s_t}+\epsilon)\}\right)\geq \epsilon \delta /4.\eqno{(3.3)}$$
We set $\{i\}=\{s_t\}$ for simplicity.
Let $I(0, 0)$ be the indicator  of the event: $\cap_i\{(l_i-\epsilon) \leq  H_{n_i}/ n_i\leq (l_i+\epsilon)\}$.
For each $\omega\in \cap_{i} \{(l_i-\epsilon) \leq  H_i/ n_i\leq (l_i+\epsilon)\}$, let $\tau_q(\omega)$ be the new configuration
by moving $\omega$ vertically  up to $q$ units. Let
$I(0, q)$ be the indicator of the  event that 
$$\{\tau_q(\omega): \omega \in \cap _i \{(l_i-\epsilon) \leq  H_{n_i}/ n_i\leq (l_i+\epsilon)\}\}.$$
We note that $\{I(0, q)\}$ is a stationary indicator sequence for $q=0, 1, \cdots$.  Since $\{t(e)\}$ is an independent and  identically distributed sequence in ${\bf E}^2$, $\{I(0, q)\}$ is an ergodic stationary sequence. By  Birkhoff's ergodic theorem and (3.3), 
$$\lim_{N\rightarrow \infty} {I(0, 0)+I(0,1)+\cdots +I(0, N) \over N} = EI(0, 0)\geq \delta \epsilon /4 \mbox{ a.s.}\eqno{(3.4)}$$
Thus, 
for any $\epsilon>0$ with $\delta > 8 \epsilon$,  by (3.4), there exists $N$ such that 
$$P(\exists \, j\mbox{ with } I(0, j)=1 \mbox { for } 1\leq j \leq N)\geq 1-\epsilon.\eqno{(3.5)}$$
For  $\omega\in \{\exists \, j\mbox{ with } I(0, j)=1 \mbox { for } 1\leq j \leq N\}$, there exists $j$, depending on $\omega$, such that there is an optimal path of $s_{0, n_i} (0, j)$ from $(0, j)$  to $\{n_i\}\times (l (1 -\epsilon)n_i+j, l (1 +\epsilon)n_i+j) $ for each $i$.

Let  ${\cal E}_{n_i}=\{\delta n_i \leq  H_{n_i}\leq n_i\}$. On ${\cal E}_{n_i}$,  there is an optimal path $\gamma_{n_i}\in \bar{S}_{n_i}$   from the origin to  ${\bf v}\in \{n_i\} \times (\delta n_i,  n_i)$ (see  Fig. 1, left).   Now we choose $n_i$ large with $l_in_i\geq  \epsilon n_i\geq 4N$  for the $l_i$ defined in (3.3) and the $N$ defined in (3.5).
Let  ${\cal E}_{n_i}(l_in_i, N)$  be the event that   there exists an optimal path $\gamma_{n_i}(l_i n_i,j)$ from $(0, l_in_i +j)$ to ${\bf o}$  for ${\bf o}\in \{n_i\} \times (-2n_i \epsilon, 2n_i\epsilon )$ and for some $0\leq j\leq N$  (see  Fig. 1, left) with passage time $s_{0, n_i}(0, l_i n_i+j)$. We want to remark that
for each $\omega\in {\cal E}_{n_i}(l_in_i,N)$, $j$ depends on $\omega$, so $(0, n_i j)$ and ${\bf o}$ may change in different configurations, but  $0\leq j\leq N$ and ${\bf o}\in \{n_i\} \times (-2n_i \epsilon, 2n_i\epsilon )$.
For each $\omega\in \{\exists \, j\mbox{ with } I(0, j)=1 \mbox { for } 1\leq j \leq N\}$ and $n_i$ selected above, we move up $\omega$ in $l_in_i$ units vertically and rotate the repositioned  configuration $180^\circ$ around the horizontal  line $y=l_in_i$  to have a new configuration in ${\cal E}_{n_i}(l_in_i,N)$.
By (3.5),  translation invariance, and symmetry,
$$ P({\cal E}_{n_i}(l_i n_i,N))\geq  P(\exists \, j\mbox{ with } I(0, j)=1 \mbox { for } 1\leq j\leq N)\geq 1-\epsilon.\eqno{(3.6)}$$
On $\{{\cal E}_{n_i} \cap {\cal E}_{n_i}(l_in_i,N)\}$, $\gamma_{n_i}$ and $\gamma_{n_i}(l_in_i,j)$ have to meet  at ${\bf w}$ for some vertex  in $L(0, n_i)$. 
Since both $\gamma_{n_i}$ and $\gamma_{n_i}(l_in_i ,j)$ are optimal paths defined above, the sub-paths from ${\bf w}$ to   ${\bf o}\in \{n_i\}\times ( -2n_i\epsilon, 2n_i\epsilon)$ along  $\gamma_{n_i}(l_in_i ,j)$ and from ${\bf w}$ to   ${\bf v}\in \{n_i\}\times (\delta n_i, n_i)$ along $\gamma_{n_i}$  have the same passage time (see Fig. 1, left). 
Thus,  on $\{{\cal E}_{n_i} \cap {\cal E}_{n_i}(l_in_i,N)\}$,  we can construct  another optimal path $\bar{\gamma}_{n_i}$  from the origin along $\gamma_{n_i}$ to meet ${\bf w}$, and along $\gamma_{n_i}(l_in_i,j)$ from ${\bf w}$ to ${\bf o}\in \{n_i\}\times (-2n_i\epsilon,  2n_i\epsilon)$ (see Fig. 1, left) with  passage time $s_{0, {n_i}}$.   We want to remark that ${\bf w}$ may change in different configurations since $j$ depends on configurations.
Let
$\bar{\cal E}_{n_i}$ be the event that there are  two optimal paths of $s_{0, n}$ from the origin to $\{n_i\}\times ( \delta n_i , 
n_i)$ and to $\{n_i\}\times (-2n_i\epsilon, 2n_i\epsilon)$, respectively. By (3.2) and (3.6),  by choosing  a small $\epsilon$ with $\delta > 8\epsilon$,  for the $n_i$ selected above,
$$P(\bar {\cal E}_{n_i})\geq P( {\cal E}_{n_i} \cap {\cal E}_{n_i}(l_i n_i,N)) \geq \delta -\epsilon\geq \delta /2.\eqno{(3.7)}$$

 \begin{figure}
\begin{center}
\setlength{\unitlength}{0.0125in}%
\begin{picture}(350,240)(67,720)
\thicklines
\put(0,770){\line(0,1){260}}
\put(0,900){\circle*{10}}
\put(0,900){\line(1,0){120}}
\put(120,900){\line(0,1){20}}
\put(120,920){\line(-1,0){30}}
\put(90,920){\line(0,1){30}}
\put(90,950){\line(1,0){90}}
\put(60,900){\line(0,-1){20}}
\put(60,880){\line(1,0){90}}
\put(150,880){\line(0,1){70}}
\put(170,875){\line(0,1){75}}
\put(170,875){\line(1,0){10}}
\put(0,970){\line(1,0){60}}
\put(60,970){\line(0,-1){90}}

\put(180,770){\line(0,1){260}}
\put(140,770){\line(0,1){260}}


\put(445,750){\mbox{$L[n_i]$}{}}

\put(180,750){\mbox{$L[n_i]$}{}}
\put(100,750){\mbox{$L[n_i(1-\kappa)]$}{}}
\put(160,950){\circle*{4}}
\put(140,950){\circle*{4}}
\put(180,950){\circle*{4}}
\put(180,875){\circle*{4}}
\put(185,873){\mbox{${\bf o}$}}
\put(190,865){\mbox{$-2n_i\epsilon$}}
\put(190,885){\mbox{$2n_i\epsilon$}}
\put(190,937){\mbox{$n_i\delta$}}
\put(190,985){\mbox{$n_i$}}
\put(-50,990){\mbox{$l_in_i+N$}}
\put(-20,945){\mbox{$l_in_i$}}
\put(0,970){\circle*{4}}
\put(-50,965){\mbox{$n_il_i+j$}}
\put(60,900){\circle*{4}}
\put(60,900){\mbox{${\bf w}$}{}}

\put(180,990){\circle{4}}
\put(180,940){\circle{4}}
\put(0,990){\circle{4}}
\put(0,950){\circle{4}}

\put(180,890){\circle{4}}
\put(180,870){\circle{4}}
\put(180,950){\circle{4}}
\put(180,950){\circle{4}}

\put(280,800){\circle*{10}}
\put(130,940){\mbox{${\bf u}$}{}}
\put(155,956){\mbox{${\bf p}$}{}}
\put(185,945){\mbox{${\bf v}$}{}}
\put(180,915){\mbox{\tiny{$(\delta -2\epsilon)n_i  \geq \delta n_i/4$}}{}}
\put(180,925){\vector(0,1){10}}
\put(180,910){\vector(0,-1){10}}


\put(20,890){\mbox{$\gamma_{n_i}$}{}}
\put(20,980){\mbox{$\gamma_{n_i}(l_in_i,j)$}{}}

\put(-20,900){\mbox{${\bf 0}$}{}}
\put(260,800){\mbox{${\bf 0}$}{}}

\put(280,800){\line(1,0){20}}
\put(300,800){\line(0,1){40}}
\put(300,840){\line(1,0){20}}
\put(320,840){\line(0,1){130}}
\put(320,970){\line(1,0){20}}
\put(340,970){\line(0,1){10}}
\put(340,980){\line(-1,0){10}}
\put(330,980){\line(0,1){20}}
\put(330,1000){\line(1,0){20}}
\put(350,1000){\line(0,-1){20}}

\put(350,980){\line(1,0){40}}
\put(390,980){\line(0,-1){170}}
\put(390,810){\line(1,0){60}}
\put(450,770){\line(0,1){260}}
\put(330,910){\dashbox(40,40)}
\put(330,910){\line(0,-1){10}}
\put(330,900){\line(1,0){80}}
\put(410,900){\line(1,0){40}}
\put(350,1000){\circle*{4}}
\put(330,970){\circle*{4}}
\put(370,980){\circle*{4}}
\put(330,910){\circle*{4}}

\put(390,900){\circle*{4}}
\put(450,800){\circle{4}}
\put(450,820){\circle{4}}
\put(450,860){\circle{4}}
\put(450,890){\circle{4}}
\put(450,910){\circle{4}}

\put(370,800){\circle{4}}
\put(350,800){\circle{4}}
\put(330,800){\circle{4}}
\put(320,956){\mbox{${\bf u}'$}{}}
\put(335,902){\mbox{${\bf z}$}{}}
\put(320,1010){\mbox{$(m_1, m_2)$}{}}
\put(347,790){\mbox{\tiny$m_1$}{}}
\put(370,790){\mbox{\tiny$m_1(1+\delta_2)$}}
\put(300,790){\mbox{\tiny$m_1(1-\delta_2)$}}
\put(370,983){\mbox{${\bf v}'$}{}}
\put(333,922){\mbox{\tiny {$B_{\delta_2 n_i}({\bf z})$}}{}}
\put(375,900){\mbox{${\bf w}$}{}}
\put(455,795){\mbox{$-\epsilon n_i$}{}}
\put(455,815){\mbox{$\epsilon n_i$}{}}
\put(455,855){\mbox{$\delta n_i/2$}{}}
\put(455,885){\mbox{$z_2-\epsilon n_i$}{}}
\put(455,905){\mbox{$z_2+\epsilon n_i$}{}}
\put(415,815){\mbox{$\gamma_{n_i}$}{}}
\put(415,905){\mbox{$\gamma({\bf z})$}{}}

\thicklines
\end{picture}
\end{center}
\caption{\em The left-hand figure is a graph proof of $H_n/n \rightarrow 0$ in probability.  If $H_n/n \not\rightarrow 0$ in probability, then there exist one optimal path $\gamma_{n_i}$ from the origin  to ${\bf v}\in \{n_i\} \times (\delta n_i, n_i)$ and another optimal path $\gamma_{n_i} (l_in_i,j)$  from $(0, n_i l_i+j)$ for $0\leq j\leq N$ to ${\bf o}\in \{n_i\}\times (-2\epsilon n_i, 2 \epsilon n_i)$ with a probability larger than $\delta /2$.  The two paths  meet at ${\bf w}$.  Both paths are optimal paths of $s_{0, n_i}(0, l_i n_i+j)$ and $s_{0, n}$.  The passage times from ${\bf w}$ along $\gamma_{n_i}$ to ${\bf v}$ and from ${\bf w}$ along $\gamma_{n_i} (l_in_i,j)$ to ${\bf o}$ are the same. There is another optimal path $\bar{\gamma}_{n_i}\in \bar{S}_{n_i}$  from the origin to ${\bf o}\in \{n_i\}\times (-2\epsilon n_i, 2\epsilon n_i)$. By using Lemma 2.3, there exists  a vertex ${\bf p}$ of a pivotal edge $e \in L(n_i(1-\kappa), n_i)$. By Lemma 2.2, $\|{\bf p}- {\bf o} \| \leq \lambda \kappa n_i$ and $\|{\bf p}- {\bf v} \| \leq \lambda \kappa n_i$.  By the triangular inequality,  $\delta n_i/4\leq \|{\bf v}-{\bf o}\|\leq  2\lambda \kappa n_i$. This is impossible if $\kappa$ is small.  $\hskip 16cm$
The right-hand figure is a graph proof of  $h^S_n/n \rightarrow 0$ in probability. 
If $h^S_n/n \not\rightarrow 0$ in probability, then there exists an optimal path $\gamma_{n_i}$ from the origin  to a vertex in
$ \{n_i\} \times (-\epsilon n_i, \epsilon n_i)$ containing $(m_1, m_2)$ with $m_2 \geq \delta n_i$ for $\delta \geq 8\epsilon$ and 
$0\leq m_1 \leq (1-\delta/(3\lambda)) n_i$ with a probability larger than $\delta/2$. By Lemma 2.2, we select $\delta_2=O(\delta)$ small such that there is a square (a dashbox)
$B_{\delta_2} n_i ({\bf z})$ below $\gamma_{n_i}$ with a probability larger than $\delta/3$.  By Lemma 3.1, for a fixed ${\bf z}$,  there is an optimal path  $\gamma({\bf z})$ from
${\bf z}=(z_1, z_2)$ to $\{n_i\}\times (z_2-\epsilon n_i, z_2+\epsilon n_i)$ such that  $z_2 -\epsilon n_i> \delta n_i /2$ with a probability larger than $1-\epsilon$. We note that $\gamma({\bf z})$ and $\gamma_{n_i}$ meet at ${\bf w}$. The passage times from ${\bf w}$ along $\gamma_{n_i}$ to $L[n_i]$ and from ${\bf w}$ along $\gamma({\bf z})$ to $L[n_i]$ are the same.
There is another optimal path in $ \bar{S}_{n_i}$  from the origin to $ \{n_i\}\times (z_2-\epsilon n_i, z_2+\epsilon n_i)$. This is impossible  if $H_i \leq \epsilon n_i$ and $z_2-\epsilon n_i \geq \delta n_i /2$ with $\delta \geq 8\epsilon$.  Since there are at most $(\lambda /\delta_2)^2$ many choices to select ${\bf z}$, the above  event  occurring with a probability is less than $(\lambda /\delta_2)^2 \epsilon$. Combining two probabilities, $\delta/3 \leq (\lambda /\delta_2)^2 \epsilon$.  The inequality cannot hold if  $\epsilon$ is small.}
\end{figure}

On $\bar{\cal E}_{n_i}$, let $\gamma_{n_i}\in \bar{S}_{n_i}$ be the optimal path from the origin to $\{n_i\} \times (\delta n_i, n_i)$, and let $\gamma({\bf u}, {\bf v})$ be the corresponding sub-path defined in Lemma 2.1.  For  the $\kappa>0$,  by (2.13) in Lemma 2.3,
\begin{eqnarray*}
&& E (D_{n_i}(\kappa); \bar{\cal E}_{n_i}) \\
&\geq & P( \tau(f)=0)  \sum_{f\in {\bf L}[0, n_i]}  P(  \exists \,\,\gamma_{n_i} \in \bar{S}_{n_i},  f\in \gamma({\bf u}, {\bf v}), \tau(f) >0, \bar{\cal E}_{n_i})\geq  F(0)   E[ N_{n_i}(\kappa);  \bar{\cal E}_{n_i}].\hskip 1cm (3.8)
\end{eqnarray*}
By (3.7), (3.8),   and (2.21),  for all large $n_i$ and $\eta$ in (2.21), there exists $C_1=C_1(F, \eta)$ such that
$$ E[D_{n_i}(\kappa); \bar{\cal E}_{n_i}]\geq F(0)   E[ N_{n_i}(\kappa);  \bar{\cal E}_{n_i}]\geq  F(0) \eta \kappa  n_i  P(N_{n_i}(\kappa)\geq \eta \kappa n_i, {\cal E}_{n_i} )\geq C_1 \kappa n_i\delta . \eqno{(3.9)}$$
By (3.9) and Lemma 2.2,  note  that $D_{n_i}(\kappa)\leq | \gamma({\bf u}, {\bf v})| $, so  there exists $C_j=C_j(F)$  for $j=2,3$ such that
\begin{eqnarray*}
C_1 \delta \kappa n_i  &\leq & E[D_{n_i}(\kappa); \bar{\cal E}_{n_i}] = \sum_{ j=1}^\infty E[D_{n_i}(\kappa); \bar{\cal E}_{n_i}\bigcap_{j=0} \{j \lambda \kappa n_i\leq D_{n_i}\leq (j+1) \lambda \kappa n_i\}]\\
&\leq &\lambda\kappa n_i P(\bar{\cal E}_{n_i}, D_{n_i}\geq 1)+\sum_{j=1}^\infty  \lambda  \kappa n_i  C_2(j+1)P( \exists\,\, \gamma_{n_i}\in \bar{S}_{n_i}: |\gamma({\bf u}, {\bf v})| \geq \lambda j \kappa n_i)\\
&\leq &\lambda\kappa n_i P(\bar{\cal E}_{n_i}, D_{n_i}\geq 1)+\sum_{j=1}^\infty  \lambda  \kappa n_i  C_2(j+1)\exp(-C_3 \kappa n_i j).\hskip 4.5cm (3.10)
\end{eqnarray*}
By (3.10),  if $n_i$ is selected to be large, then there exists $C_4 =C_4(F, \delta)>0$ independent  of $\kappa$ and $n_i$ such that
$$ C_4 \leq P(\bar{\cal E}(n_i)\cap \{\exists \mbox{  a pivotal edge $e$ of $S_{n_i}$ in } L(n_i-\kappa n_i, n_i)\}).\eqno{(3.11)}$$
On $ \{\exists \mbox{  a pivotal edge }e\in L(n_i-\kappa n_i, n_i)\}$,  we assume that ${\bf p}\in \gamma({\bf u}, {\bf v})$ is a vertex of $e$ and every optimal path in $\bar{S}_{n_i}$ has to pass through ${\bf p}$ (see  Fig. 1, left).
By Lemma 2.2,  by taking $n_i$ large, there exists $ C_5=C_5(F)$ such that
$$ P( d({\bf p}, {\bf v}) \leq  \lambda \kappa n_i  \mbox{ and } d({\bf p}, {\bf o}) \leq \lambda \kappa n_i )\geq C_5/4.\eqno{(3.12)}$$
By  the triangular inequality, with a positive probability, 
$$\delta n_i /4\leq n_i\delta - 2n_i \epsilon =n_i (\delta -2\epsilon) < 2\lambda \kappa n_i.\eqno{(3.13)}$$
Since $\lambda$ and $\delta$ are assumed to be  positive constants independent of $\kappa$,  (3.13) will contradict  if $\kappa$  is small.
The contradiction tells us that (3.1) cannot hold. Therefore, Lemma 3.1 follows. $\Box$\\ 

Let 
$$h_n^S=\max\{d({\bf u}, \mbox{ the $X$-axis}) : {\bf u}\in S_n\}.$$
We show the following lemma.\\

{\bf Lemma 3.2.} {\em If  $F$ is a  distribution with $0< F(0) < p_c$ and (1.5) holds, then
$$\lim_{n\rightarrow \infty}  {{h}_n^S\over n}=0\mbox{ in probability}.\eqno{}$$}

{\bf Proof.}  
If we suppose that Lemma 3.2 does not hold, then by symmetry, then there exists a subsequence $\{n_i\}$ such that
$$P(  h^S_{n_i} /n_i \geq \delta) \geq \delta \mbox{ for all large }\{n_i\}.\eqno{(3.14)}$$
On $h^S_{n_i} /n_i \geq \delta$ for all $i$, there exist $\gamma_{n_i} \in \bar{S}_{n_i}$ and   $(m_1, m_2)\in \gamma_{n_i}$  such that
$$m_2 =h^S_{n_i}\geq \delta n_i \mbox{ for  all } i.$$
  If there are many $\{(m_1, m_2)\}$  such that $h^S_{n_i} =m_2 \geq \delta n_i$,  then we simply select  $(m_1, m_2)$ among all the optimal paths 
  $\{\gamma_{n_i}\}$ in $\bar{S}_{n_i}$ with the  largest $x$-coordinate  when we go along optimal paths starting from the origin. We also select $\gamma_{n_i}$ containing $(m_1, m_2)$.
Thus,  by Lemma 3.1 and the assumption of (3.14), for  a small $\epsilon$ with $\delta > 8\epsilon >0$,
$$P(\exists \,\,\, (m_1, m_2) \mbox{ such that } m_2 =h^S_{n_i} \geq \delta n_i, H_{n_i}/ n_i\leq \epsilon ) \geq \delta/2 \mbox{ for all  large }\{n_i\}.\eqno{(3.15)}$$
On $\{\exists \,\,\, (m_1, m_2) \mbox{ such that } m_2 =h^S_{n_i} \geq \delta n_i, H_{n_i}/ n_i\leq \epsilon\}$,  we know that $0< m_1< n_i$.
We assume that  $n_i-m_1 = \delta_1 n_i$ for $\delta_1>0$.  We will show that $\delta_1=O(\delta)$ by using Lemma 2.2. 
By Lemma 2.2,  with a probability larger than $1-\exp(-O(\delta n_i))$, ${\bf w}=(w_1, w_2)\in L[n_i]\cap \gamma_{n_i}$,  the terminate vertex of $\gamma_n$, such that
$$ w_2\leq H_{n_i} \mbox{ and } d((m_1, m_2), {\bf w}) \leq \lambda  \delta_1 n_i.\eqno{(3.16)}$$
Note that $m_2 > \delta n_i$ and $H_{n_i}\leq \epsilon n_i$, so 
$$\delta  n_i/2 \leq d((m_1, m_2), {\bf w})\leq \lambda \delta_1 n_i.\eqno{(3.17)}$$
By (3.17), 
$$P(\delta_1 \geq \delta / (2\lambda))\geq 1-\exp(-O(\delta n_i)).\eqno{(3.18)}$$
 Similarly, we show that  
 $$P( n_i \delta /(2\lambda)\leq m_1\leq  n_i(1-\delta /(2\lambda)))\geq 1-\exp(-O(\delta n_i)).\eqno{(3.19)}$$
By (3.19), we choose $\delta_2=O(\delta)$  small,  independent of $n_i$, such that 
$$P( n_i \delta /(3\lambda) \leq  m_1-\delta_2 n_i\leq m_1 \leq    m_1+\delta_2 n_i \leq n_i(1-\delta/(3\lambda)))\geq 1-\exp(-O(\delta n_i)).\eqno{(3.20)}$$
Let  $\gamma_{n_i}$, starting from the origin,  first meet  $L[m_1-\delta_2 n_i]$ at ${\bf u}'$ and last meet $L[m_1+n_i\delta_2]$ at ${\bf v}'$ (see  Fig. 1, right).
Let $\gamma({\bf u}', {\bf v}')$ be the sub-path of $\gamma_{n_i}$ from ${\bf u}'$ and ${\bf v}'$ as we defined in Lemma 2.1.
   By   Lemma 2.2 again, for all large $n_i$, 
  $$P(|\gamma({\bf u}', {\bf v}')|\leq \lambda \delta_2 n_i) \geq 1-\exp(-O(\delta_2 n_i)).\eqno{(3.21)}$$
  By (3.15), (3.20), and (3.21), for large $n_i$, we select $\delta_2=O(\delta)$  small independent of $n_i$ such that
  \begin{eqnarray*}
  && P(\exists \,\,\, (m_1, m_2) \mbox{ such that } m_2 =h^S_{n_i} \geq \delta n_i, 0< m_1-\delta_2n_i\leq   m_1+\delta_2 n_i\leq n_i(1-\delta/(3\lambda)),\\
  && \hskip 8cm  |\gamma({\bf u}', {\bf v}')|\leq \lambda \delta_2 n_i, H_{n_i}/ n_i\leq \epsilon ) \geq \delta/3.\hskip 1.3cm {(3.22)}
  \end{eqnarray*}
Note that $(m_1, m_2) \in \gamma({\bf u}', {\bf v}')$, so if $\delta_2$ is selected to be small but independent  of $n_i$,   by (3.22), 
\begin{eqnarray*}
&&P(\exists \,\,\gamma_{n_i} \in \bar{S}_{n_i}, (m_1, m_2) \in \gamma_{n_i}  \mbox{ such that }\forall (x_1, x_2)\in \gamma({\bf u}', {\bf v}') \\
&&\hskip 2cm \mbox{ with  } 0\leq x_1 \leq n_i(1-\delta/(4\lambda)), x_2 > n_i\delta /2, H_{n_i}/n_i \leq \epsilon ) \geq \delta /3.\hskip 3cm {(3.23)}
\end{eqnarray*}

Now we divide ${\bf Z}^2$ into equal disjoint squares,  except their boundaries, with length $\delta_2 n$ for small $\delta_2 < \delta/4$  with  $\delta_2=O(\delta)$.  More precisely,  for ${\bf z}=(z_1, z_2)\in {\bf Z}^2$,
$$B_{\delta_2 n} ({\bf z})=[ \delta_2 n z_1, 2\delta_2 n z_1]\times [ \delta_2 n z_2, 2\delta_2 n z_2].$$
If $\delta_2$ is small, then on $\{\forall (x_1, x_2)\in \gamma({\bf u}, {\bf v}), 0\leq  x_1 \leq n_i(1-\delta/(4\lambda),  x_2 > n_i\delta /2\}$,
there exists a square $B_{\delta_2 n_i} ({\bf z})$ {\em below} $\gamma_{n_i}$  with ${\bf z}=(z_1, z_2)$ such that  (see Fig. 1, right)
$$0\leq z_1 \leq n_i(1-\delta/(4\lambda))\mbox{ and } z_2 \geq \delta n_i/4.\eqno{(3.24)}$$
A square below a curve means that any vertical line from $-\infty$ to $\infty$ always meets the square first before the curve.
For a ${\bf z}$ satisfying (3.24), let ${\cal E}_{n_i}(\delta, {\bf z})$ be the event of  all optimal paths $\gamma({\bf z})$ from ${\bf z}=(z_1, z_2) $ to $L[n_i]$ with $z_1 < n_i(1-\delta/(4\lambda))$ and $z_2 \geq \delta n_i/4$,  and with a passage time 
$s_{z_1, n_i}({\bf z})$  from ${\bf z}$ to $\{n_i\}\times [z_2-\epsilon n_i, z_2+\epsilon n_i]$ in $L(z_1, n_i)$ (see Fig. 1, right). By Lemma 3.1,  for any small $\epsilon >0$ and for a fixed ${\bf z}$, if  $n_i$ is large, then 
$$P({\cal E}_{n_i}(\delta, {\bf z}))\geq 1-\epsilon.\eqno{(3.25)}$$
Now we need to choose   $B_{\delta_2 n_i} ({\bf z})$ from $\{B_{\delta_2 n_i} ({\bf u})\}$.
On $\gamma_n\subset [0, n_i]\times [-\lambda n_i, \lambda n_i]$, there are at most $(\lambda /\delta_2)^2$ many choices for $\{B_{\delta_2 n_i} ({\bf z})\}$.  When $B_{\delta_2 n_i} ({\bf z})$ is fixed, so is ${\bf z}$. On ${\cal E}_{n_i}(\delta_2, {\bf z})\cap \{H_{n_i} \leq \epsilon n_i\}$, $\gamma({\bf z})$ and $\gamma_{n_i}$  have to meet ${\bf w}$ in
$L(0, n_i)$ (see Fig. 1, right). So we can construct another optimal path from the origin along $\gamma_{n_i}$ to ${\bf w}$  and along $\gamma({\bf z})$ from ${\bf w}$
to $\{n_i\}\times [z_2-\epsilon n_i, z_2+\epsilon n_i]$ with passage time $s_{0, n}$. 
Since $z_2 \geq n_i \delta /4$,  
$$H_{n_i} \geq z_2-\epsilon n_i\geq 2\epsilon n_i.\eqno{(3.26)}$$
This will contradict  that $H_{n_i}\leq  \epsilon  n_i$ (Fig. 1, right). So  for a fixed ${\bf z}$,
$${\cal E}_{n_i}(\delta, {\bf z})\cap\{ H_{n_i} \leq \epsilon n_i\}=\emptyset.\eqno{(3.27)}$$
By (3.18), (3.25), and (3.27), for $\delta_2$, $\delta$, and $\epsilon$ defined above,
\begin{eqnarray*}
&\delta /3 &\leq P(\exists \,\,\gamma_{n_i} \mbox{ with } T(\gamma_{n_i})= s_{0, {n_i}}, (m_1, m_2) \in \gamma_{n_i}  \mbox{ such that }\forall (x_1, x_2)\in \gamma({\bf u}, {\bf v}) \\
&&\hskip 4cm \mbox{ with  }  x_1 \leq n_i(1-\delta/(4\lambda)), x_2 > n_i\delta /2, H_{n_i}/n_i \leq \epsilon )\\
&\leq & (\lambda /\delta_2)^2 [P(\exists \,\,\gamma_{n_i} \mbox{ with } T(\gamma_{n_i})= s_{0, n_i}, (m_1, m_2) \in \gamma_{n_i}  \mbox{ such that }\forall (x_1, x_2)\in \gamma({\bf u}, {\bf v})\\
&&\hskip 1cm \mbox{ with  }  x_1 \leq n_i(1-\delta/(4\lambda)), x_2 > n_i\delta /2, {\cal E}_{n_i}(\delta, {\bf z})\mbox{ for a fixed }{\bf z}, H_{n_i}/n_i \leq \epsilon)+\epsilon]\\
&\leq &  (\lambda /\delta_2)^2 \epsilon.\hskip 14cm (3.28)
\end{eqnarray*}
Since $\delta$, $\delta_2$, and $\lambda$ are  independent from $n_i$, (3.28) cannot hold if $\epsilon$ is small for a large $n_i$. The contradiction tells that  (3.14) cannot hold. Lemma 3.2 follows. $\Box$\\

{\bf Proof of theorem.}   
 \begin{figure}
\begin{center}
\setlength{\unitlength}{0.0125in}%
\begin{picture}(250,240)(67,720)
\thicklines
\put(117,770){\line(0,1){260}}
\put(350,770){\line(0,1){260}}
\put(234,770){\line(0,1){260}}
\put(0,770){\line(0,1){260}}
\put(117,900){\circle*{10}}
\put(105,885){{\bf 0}}
\put(0,900){\dashbox{1}(350, 0.5)}
\put(0,850){\dashbox{1}(350, 0.5)}
\put(0,950){\dashbox{1}(350, 0.5)}
\put(117,900){\line(-1,0){20}}
\put(97,900){\line(0,1){10}}
\put(97,910){\line(1,0){60}}
\put(157,910){\line(0,1){60}}
\put(157,970){\line(1,0){30}}
\put(187,970){\line(0,-1){100}}
\put(187,870){\line(1,0){60}}
\put(247,870){\line(0,1){30}}
\put(247,900){\line(-1,0){10}}

\put(0,925){\line(1,0){130}}
\put(130,925){\line(0,1){20}}
\put(130,945){\line(1,0){220}}

\put(0,875){\line(1,0){150}}
\put(150,875){\line(0,-1){10}}
\put(150,865){\line(1,0){200}}

\put(230,750){\mbox{$L[n]$}{}}
\put(110,750){\mbox{$L[0]$}{}}
\put(330,750){\mbox{$L[\lambda n]$}{}}
\put(-10,750){\mbox{$L[-\lambda n]$}{}}
\put(158,944){\circle*{4}}
\put(188,944){\circle*{4}}

\put(-25,950){\mbox{$2\epsilon  n$}}
\put(-35,875){\mbox{$-\epsilon  n$}}
\put(-20,920){\mbox{$\epsilon  n$}}
\put(-35,850){\mbox{$-2\epsilon  n$}}
\put(0,925){\circle*{10}}
\put(0,875){\circle*{10}}

\put(234,900){\circle*{10}}
\put(234,905){\mbox{$(n,0)$}}

\put(25,930){\mbox{$\gamma_{3\lambda n}(\epsilon n)$}}
\put(25,880){\mbox{$\gamma_{3\lambda n}(-\epsilon n)$}}
\put(130,913){\mbox{$\gamma_n$}}
\put(145,935){\mbox{${\bf u}$}}
\put(190,935){\mbox{${\bf v}$}}

\thicklines
\end{picture}
\end{center}
\caption{\em  The figure, a graph proof of the theorem,  shows that $h_n^A\geq 2\epsilon$ has a small probability. ${\cal E}_n^+$  is the event that  all the optimal paths 
$\{\gamma_{3\lambda}(\epsilon n)\}$ of $s_{-\lambda n, \lambda n} (-\lambda n, \epsilon n)$ stay inside $[-\lambda n, \lambda n]\times [0, 2\epsilon n]$. ${\cal E}_n^-$  is the event that  all the optimal paths $\{\gamma_{3\lambda}(-\epsilon n)\}$ of $s_{-\lambda n, \lambda n} (-\lambda n, -\epsilon n)$ stay inside $[-\lambda n, \lambda n]\times [-2\epsilon n,0]$. By Lemma 3.2, 
${\cal E}^+_n\cap {\cal E}^-_n$ has  a probability  larger than $1-\epsilon$ for all large $n$.   With a probability larger than $1-\exp(-O(n))$, all optimal paths
of $a_{0, n}$ stay inside $[-\lambda n, \lambda n]^2$, denoted by ${\cal E}_n$, the event.
If  $\{h_n^A\geq 2\epsilon\}\cap {\cal E}_n$, then there is an optimal path $\gamma_n$ 
crossing   the horizontal line  $y=2\epsilon n$ or $y=-2\epsilon n$.  On ${\cal E}^+_n\cap {\cal E}^-_n$, $\gamma_n$ has to meet an optimal path
$\gamma_{3\lambda n}(\epsilon n)$ or $\gamma_{3\lambda n}(-\epsilon n)$. In the first case, $\gamma_n$ and $\gamma_{3\lambda n}(\epsilon n)$
meet at ${\bf u}$ and ${\bf v}$. Then we can construct an optimal path such that it goes along $\gamma_{3\lambda n}(\epsilon n)$ to meet ${\bf u}$,  then along $\gamma_n$ to meet ${\bf v}$,  and finally along $\gamma_{3\lambda n}(\epsilon n)$ from ${\bf v}$ to meet $L[\lambda n]$. The constructed optimal path will not stay inside $[-\lambda n, \lambda n]\times [0, 2\epsilon n]$.  So ${\cal E}^+_n$ cannot occur. $h_n^A\geq 2\epsilon$ has a  probability less than $\epsilon$. }
\end{figure}
For the $\lambda$ in (1.8) and (1.9), let
$${\cal E}_n=\{ \forall \gamma_n \subset L(-\lambda n, \lambda n) \mbox{ with }T(\gamma_n)= a_{0, n} \mbox{ or } b_{0, n})\}.$$
By (1.8) and (1.9) with the   $\lambda$,
 $$P({\cal E}_n)\geq 1-\exp(-O(n)).\eqno{(3.29)}$$
 We denote by  (see Fig. 2)
 ${\cal E}_n^+(\epsilon)$  the event that optimal paths of $s_{-\lambda n, n}(-\lambda n, \epsilon n)$ stay inside $[-\lambda n, n]\times [0, 2\epsilon n]$.
 Moreover,  we also denote by (see Fig. 2) ${\cal E}_n^-(\epsilon)$  the event that optimal paths of $s_{-\lambda n, n}(-\lambda n, -\epsilon n)$ stay inside $[-\lambda n, n]\times [-2\epsilon n,0]$. If $n$ is large, then by  Lemma 3.2, 
 $$P({\cal E}^+_n(\epsilon)\cap {\cal E}^-_n(\epsilon))\geq 1-\epsilon.\eqno{(3.30)}$$
 On ${\cal E}_n \cap {\cal E}^+_n(\epsilon)\cap {\cal E}^-_n(\epsilon)$, by the graph proof in Fig. 2,  optimal paths of $a_{0, n}$ or $b_{0, n}$ have to stay inside $[-\lambda n, n]\times [-2\epsilon n,2\epsilon n]$, otherwise either ${\cal E}^+_n(\epsilon)$ or ${\cal E}^-_n(\epsilon)$ cannot occur.
 By (3.29) and (3.30), for large $n$, 
 $$ P(|h^{\alpha}_n|\geq 2\epsilon n) \leq P(|h^{\alpha}_n|\geq 2\epsilon n, {\cal E}_n \cap {\cal E}^+_n(\epsilon)\cap {\cal E}^-_n(\epsilon))+P(({\cal E}_n \cap {\cal E}^+_n(\epsilon)\cap {\cal E}^-_n(\epsilon))^C)\leq  2\epsilon.\eqno{(3.31)}$$
Thus, by (3.31), $h_n^\alpha/n$  converges to zero  in probability in the theorem for $\alpha=A$ or $B$.  

It remains to show that $h_n^\alpha/n$  converges to zero  in $L_1$. By (3.31), (1.8) and (1.9),   for $\epsilon >0$, 
$$E|h_n^\alpha|/ n=  E( |h_n^\alpha|/ n; |h_n^\alpha |\leq \lambda n))+ \sum_{m\geq \lambda n} (\lambda m)^2 \exp(-O(m)) \leq 2\epsilon + 2\epsilon \lambda +  \exp(-O(n)).\eqno{(3.32)}$$
By (3.32), $h_n^\alpha/n$ converges to zero in $L_1$  for $\alpha=A$ or $B$. 
$\Box$

\newpage

\begin{center}
{\bf \large References}
\end{center}
  Berg, J. van den and Kesten, H. (1993).  Inequalities for the time constant in first-passage percolation. {\em Ann. Appl. Probab.} {\bf 3}
56--80.\\
Chow, Y. and Zhang, Y. (2003). Large deviation in first passage percolation. {\em Ann. Appl. Probab.} {\bf 4}
1601--1614.\\
Cox, J. T.  and Durrett, R. (1981). Some limit theorems for percolation processes with necessary and sufficient condition. {\em Ann. Probab.} {\bf 9} 583--603.\\
Durrett, R. and  Liggett, T. (1981).
The shape of the limit set in Richardson's growth model. {\em Ann.
Probab.} {\bf 9} 186--193.\\
Grimmett, G. (1999). {\em Percolation.} Springer, Berlin.\\
 Hammersley, J. M. and  Welsh, D. J. A. (1965).
First-passage percolation, subadditive processes,
stochastic networks and generalized renewal theory.
In {\em Bernoulli, Bayes, Laplace Anniversary Volume} 
(J. Neyman   and L. LeCam,  eds.) 61--110. Springer, Berlin.\\
Kesten, H. (1986). Aspects of first-passage percolation. {\em Lecture Notes in
Math.} {\bf 1180} 125--264. Springer, Berlin.\\
Kesten, H. (1993). On the speed of convergence in first passage percolation. {\em Ann.  Appl. Probab.} {\bf 3} 296--338.\\
Nakajima, S. (2019).  On properties of optimal paths in first-passage perccolation. {\em J.   Stat.  Phys.} {\bf 174} 259--275.\\
Smythe, R. T. and Wierman, J. C. (1978).
First passage percolation on the square lattice.
{\em Lecture Notes in Math.} {\bf 671} Springer, Berlin.\\
 Zhang, Y. (2006). The divergence of fluctuations for shape in first passage percolation. {\em Probab. Theory Relat.  Fields}.
{\bf 136} 298--320.\\
Zhang, Y. (2010).   On the concentration and the convergence rate with a moment condition in first passage percolation.
{\em  Stochastic Process Appl.} {\bf 120} 1317--1341.\\

\noindent
Yu Zhang\\
Department of Mathematics\\
University of Colorado,
Colorado Springs, CO 80933\\
email: yzhang3@uccs.edu\\
\end{document}